\documentclass{article}
\usepackage{amssymb,amsmath,bm,geometry,graphics,graphicx,color,url,hyperref}

\def\no{\noindent}
\def\pmatrix{\left(\begin{array}}
\def\endpmatrix{\end{array}\right)}

\def\RR{\mathbb{R}}

\def\I{{\cal I}}

\def\P{{\cal P}}

\def\dd{\mathrm{d}}

\newtheorem{theo}{Theorem}

\newtheorem{rem}{Remark}

\def\proof{\noindent\underline{Proof}\quad}
\def\QED{\mbox{~$\Box{~}$}}

\def\bfb{{\bm{b}}}
\def\bfc{{\bm{c}}}

\def\bfgamma{{\bm{\gamma}}}

\begin{document}
\title{A note on the continuous-stage Runge-Kutta(-Nystr\"om) formulation of Hamiltonian Boundary Value Methods (HBVMs)}

\author{Pierluigi~Amodio$^1$ \and Luigi~Brugnano$^2$ \and Felice~Iavernaro$^1$}
\date{$^1$\, {\scriptsize Dipartimento di Matematica, Universit\`a di Bari, Italy, \url{{pierluigi.amodio,felice.iavernaro}@uniba.it} }\\
$^2$\, {\scriptsize Dipartimento di Matematica e Informatica ``U.\,Dini'', Universit\`a di Firenze, Italy, \url{luigi.brugnano@unifi.it}}}
\maketitle

\begin{abstract} In recent years, the class of energy-conserving methods named Hamiltonian Boundary Value Methods (HBVMs) has been devised for numerically solving Hamiltonian problems. In this short note, we study their natural formulation as continuous-stage Runge-Kutta(-Nystr\"om) methods, which allows a deeper insight in the methods.
\medskip

\no{\bf Keywords:} continuous-stage Runge-Kutta methods,  Runge-Kutta-Nystr\"om methods, Hamiltonian Boundary Value Methods, HBVMs.
\medskip

\no{\bf MSC:} 65L05, 65P10.

\end{abstract}

\section{Introduction}

The numerical solution of Hamiltonian problems has been recently tackled by defining energy-conserving methods, which can be regarded as continuous-stage Runge-Kutta (RK, hereafter) methods (e.g., \cite{QMcL2008,BIT2010,Ha2010,TS2014}). In their simplest (and most effective) form,\footnote{For a more general form which, however, we shall not consider here,  we refer, e.g., to \cite{MiBu2016}.}   continuous-stage RK methods are ``methods'' that, when applied for solving an initial value problem for ODEs (ODE-IVP, hereafter), which we assume without loss of generality in the form
\begin{equation}\label{ivp}
\dot y(t) = f(y(t)), \qquad t\in[0,h], \qquad y(0) = y_0\in\RR^m,
\end{equation}
with $f$ analytical, define an approximating function $u:[0,h]\rightarrow \RR^m$ such that
\begin{equation}\label{csRK1}
u(ch) = y_0 + h\int_0^1 a_{c\tau} f(u(\tau h))\dd\tau, \qquad c\in[0,1],
\end{equation}
with $a_{c\tau}:[0,1]\times[0,1]\rightarrow\RR$, and a corresponding approximation to $y(h)$,
\begin{equation}\label{csRK2}
y_1 = y_0 + h\int_0^1 f(u(ch))\dd c.
\end{equation} 
As is usual, this procedure can be summarized by the following (generalized) Butcher tableau,
$$\begin{array}{c|c}
c &a_{c\tau}\\[2mm]
\hline
\\[-2mm]
& 1\end{array}~.
$$
We observe that (\ref{csRK1})-(\ref{csRK2}) is not yet an actual numerical method, due to the fact that the involved  integrals  need to be conveniently approximated by means of quadrature rules. In so doing, one obtains ``usual'' RK methods.\footnote{I.e., having discrete stages.} Nevertheless, (\ref{csRK1})-(\ref{csRK2}) can be useful for purposes of analysis \cite{Bu2008,BIT2012,TS2012,TS2014,TZ2018,TSZ2019,LiWu2019} since, essentially, it allows to discuss all Runge-Kutta methods derived by using different quadratures for approximating the involved integrals.  In particular, the papers \cite{TZ2018,TSZ2019} have inspired the present note, where we provide the continuous-stage RK formulation of Hamiltonian Boundary Value Methods (HBVMs) \cite{BIT2009,BIS2010,BIT2010,BIT2012-1,BIT2015,LIMbook2016,BI2018}, a class of energy-conserving methods for Hamiltonian problems, which have been developed along several directions \cite{BCMR2012,BI2012,BIT2012-2,BIT2012-3,BS2014,BGI2018,BGIW2018}, including Hamitonian BVPs \cite{ABI2015}, highly-oscillatory problems \cite{BMR2018,ABI2019}, Hamiltonian PDEs \cite{BFCI2015,BBFCI2018,BIMR2018,BFCI2019,BGS2019,BZL2018,BGZ2019}, and also considering their efficient implementation \cite{BIT2011,BFCI2014}. Here, we shall also consider the continuous formulation of such methods when applied for solving special second-order problems  \cite{LIMbook2016}, i.e., problems in the form
\begin{equation}\label{ivp2}
\ddot q(t) = f(q(t)), \qquad t\in[0,h], \qquad q(0) = q_0, ~ \dot q(0) = p_0\,\in\,\RR^m,
\end{equation}
where, for the sake of brevity,  we shall again assume $f$ to be analytical.

With these premises, the structure of the paper is as follows: in Section~\ref{fop} we study the case of first order ODE problems; Section~\ref{sop} is devoted to study the case where one solves special second-order problems; at last, a few concluding remarks are drawn in Section~\ref{fine}.

\section{The framework}\label{fop}

Generalizing the arguments in \cite{BIT2012}, let us consider the orthonormal Legendre polynomial basis $\{P_j\}_{j\ge0}$ on the interval $[0,1]$:
\begin{equation}\label{orto}
P_j\in\Pi_j, \qquad \int_0^1 P_i(c)P_j(c)\dd c = \delta_{ij}, \qquad \forall i,j\ge0,
\end{equation} where $\Pi_j$ is the set of polynomials of degree $j$. 
Then, the ODE-IVP (\ref{ivp}) can be written, by expanding the right-hand side along the Legendre basis, as
\begin{equation}\label{ivp1}
\dot y(ch) = \sum_{j\ge0} P_j(c)\gamma_j(y), \qquad c\in[0,1], \qquad \gamma_j(y) = \int_0^1P_j(\tau) f(y(\tau h))\dd \tau, \qquad j\ge0,
\end{equation}
from which, integrating side by side, one obtains the following formal expression for the solution of (\ref{ivp}):
\begin{equation}\label{y}
y(ch) = y_0 + h\sum_{j\ge0} \int_0^c P_j(x)\dd x\, \gamma_j(y), \qquad c\in[0,1].
\end{equation}
The above equations can be cast in vector form by introducing the infinite vectors
\begin{equation}\label{PIgamma}
\P_\infty(c) := \pmatrix{c} P_0(c)\\[1mm] P_1(c)\\ \vdots\, \endpmatrix,\quad 
 \I_\infty(c) := \int_0^c \P_\infty(x)\dd x \equiv \pmatrix{c} \int_0^cP_0(x)\dd x\\[2mm]  \int_0^cP_1(x)\dd x\\ \vdots\, \endpmatrix,\quad
 \bfgamma(y) = \pmatrix{c} \gamma_0(y)\\[1mm] \gamma_1(y)\\ \vdots \endpmatrix,
\end{equation}
respectively as:
\begin{eqnarray}\label{y1gamma}
\dot y(ch) = \P_\infty(c)^\top\otimes I_m\bfgamma(y), \qquad c\in[0,1],\qquad  \bfgamma(y) = \int_0^1 \P_\infty(\tau)\otimes I_m f(y(\tau h))\dd \tau,
\end{eqnarray}
and
\begin{equation}\label{ych}
y(ch) = y_0 + h\I_\infty(c)^\top\otimes I_m\,\bfgamma(y), \qquad c\in[0,1].
\end{equation}
Moreover, by considering that
\begin{equation}\label{IXPPI}
\I_\infty(c)^\top = \P_\infty(c)^\top X_\infty, \qquad \int_0^1 \P_\infty(c) \P_\infty(c)^\top\dd c = I,
\end{equation}
with $I$ the identity operator and
\begin{equation}\label{X}
X_\infty = \pmatrix{cccc} 
\xi_0 & -\xi_1\\
\xi_1 & 0   & -\xi_2 \\
         & \xi_2 & \ddots &\ddots\\
         &          &\ddots  &\ddots\endpmatrix, \qquad \xi_i = \frac{1}{2\sqrt{|4i^2-1|}}, \quad i\ge0,
         \end{equation}
one also obtains that
\begin{equation}\label{PIX}
\int_0^1 \P_\infty(c) \I_\infty(c)^\top\dd c = X_\infty.
\end{equation}

Setting \,$y_1\equiv y(h),$\, we can cast (\ref{ych}) as:
\begin{eqnarray}\nonumber
y(ch) &=& y_0 + h\int_0^1 \I_\infty(c)^\top \P_\infty(\tau)f(y(\tau h))\dd\tau, \qquad c\in[0,1],\\
y_1   &=& y_0 + h\int_0^1 f(y(ch))\dd c,\label{cRK}
\end{eqnarray}
which, by virtue of (\ref{IXPPI})--(\ref{PIX}), can be also written as
\begin{eqnarray}\nonumber
y(ch) &=& y_0 + h\int_0^1 \P_\infty(c)^\top X_\infty \P_\infty(\tau)f(y(\tau h))\dd\tau, \qquad c\in[0,1],\\
y_1   &=& y_0 + h\int_0^1 f(y(ch))\dd c.\label{cRK1}
\end{eqnarray}
In other words, we are speaking about the application of the following continuous-stage RK method to problem (\ref{ivp}) :
\begin{equation}\label{RK}
\begin{array}{c|c}
c & \I_\infty(c)^\top \P_\infty(\tau)\\[2mm]
\hline
\\[-2mm]
& 1 \end{array}
\qquad \equiv \qquad
\begin{array}{c|c}
c & \P_\infty(c)^\top X_\infty \P_\infty(\tau)\\[2mm]
\hline
\\[-2mm]
& 1\end{array}\qquad =: \qquad
\begin{array}{c|c}
c &a_{c\tau}^{(\infty)}\\[2mm]
\hline
\\[-2mm]
& 1\end{array}~.
\end{equation}
As is clear,  by virtue of (\ref{IXPPI})-(\ref{X}), the coefficients of this ``continuous-stage RK method'', providing the exact solution of (\ref{ivp}), are given by
\begin{eqnarray}\nonumber
a_{c\tau}^{(\infty)} &:=& \sum_{j=0}^\infty \int_0^c P_j(x)\dd x\, P_j(\tau) \\ \label{actau}
&\equiv& c + \sum_{j=1}^\infty [\xi_{j+1}P_{j+1}(c)-\xi_jP_{j-1}(c)]P_j(\tau), \qquad c,\tau\in[0,1].
\end{eqnarray}

\subsection{Polynomial approximation}\label{poli1}
In order to obtain a polynomial approximation $\sigma\in\Pi_s$ to $y$, let us now introduce the truncated vectors 
\begin{equation}\label{PIs} 
\P_s(c) := \pmatrix{c} P_0(c) \\ \vdots\\ P_{s-1}(c)\endpmatrix, \qquad
 \I_s(c) := \int_0^c \P_s(x)\dd x \equiv \pmatrix{c} \int_0^cP_0(x)\dd x \\ \vdots\\ \int_0^cP_{s-1}(x)\dd x\endpmatrix,
\end{equation}
in place of the corresponding infinite ones in (\ref{PIgamma}). In so doing, we replace (\ref{RK}) with the continuous-stage RK method 
\begin{equation}\label{RKs}
\begin{array}{c|c}
c & \I_s(c)^\top \P_s(\tau)\\[2mm]
\hline
\\[-2mm]
& 1 \end{array}
\qquad  =: \qquad
\begin{array}{c|c}
c &a_{c\tau}^{(s)}\\[2mm]
\hline
\\[-2mm]
& 1\end{array}~,
\end{equation}
whose coefficients are now polynomials of degree $s$. Consequently, by setting now $y_1\equiv \sigma(h)$ the approximation to $y(h)$, one obtains:
\begin{equation}\label{RKs1}
\sigma(ch) = y_0+h\int_0^1 a_{c\tau}^{(s)}f(\sigma(\tau h))\dd\tau, \qquad c\in[0,1], \qquad y_1 = y_0+h\int_0^1 f(\sigma(ch))\dd c.
\end{equation}
The following straightforward result holds true.

\begin{theo}\label{HBVMs} The continuous-stage RK method (\ref{RKs})-(\ref{RKs1}) coincides with the HBVM$(\infty,s)$ method in \cite{BIT2010}.\footnote{In particular when $s=1$ one retrieves the AVF method in \cite{QMcL2008}.}\end{theo}
\proof In fact, from (\ref{PIs}), one has that (\ref{RKs1}) is equivalent to
\begin{equation}\label{MFE}
\sigma(ch) = y_0 + h\sum_{j=0}^{s-1} \int_0^cP_j(x)\dd x \int_0^1 P_j(\tau)f(\sigma(\tau h))\dd\tau, \qquad c\in[0,1],
\end{equation}
which, according to \cite[Definition\,1]{BIT2010} (see also \cite{BI2018}), is the {\em Master Functional Equation} defining a HBVM$(\infty,s)$ method.\,\QED\medskip

Furthermore, by considering that (see (\ref{PIs}) and (\ref{X}))
\begin{equation}\label{PIXs}
\I_s(c)^\top = \P_{s+1}(c)^\top \pmatrix{ccccc} 
\xi_0 & -\xi_1\\
\xi_1 & 0   & -\xi_2 \\
         & \xi_2 & \ddots &\ddots\\
         &          &\ddots  &\ddots &-\xi_{s-1}\\
         &          &            & \xi_{s-1} & 0\\ \hline 
         &          &            &               &\xi_s \endpmatrix =:\P_{s+1}(c)^\top\hat X_s \equiv \P_{s+1}(c)^\top \pmatrix{c} X_s\\ \hline 0,\dots,0,\xi_s \endpmatrix,
\end{equation}
one easily obtains that (compare with (\ref{actau}))
\begin{eqnarray}\label{actaus}
a_{c\tau}^{(s)} &=& \sum_{j=0}^{s-1} \int_0^c P_j(x)\dd x\, P_j(\tau) \\ \nonumber
&\equiv& c + \sum_{j=1}^{s-1} [\xi_{j+1}P_{j+1}(c)-\xi_jP_{j-1}(c)]P_j(\tau), \qquad c,\tau\in[0,1].
\end{eqnarray}
As a result, from (\ref{RKs}) and (\ref{PIXs}), one obtains that
\begin{equation}\label{RK1s}
\begin{array}{c|c}
c & \I_s(c)^\top \P_s(\tau)\\[2mm]
\hline
\\[-2mm]
& 1 \end{array}
\qquad \equiv \qquad
\begin{array}{c|c}
c & \P_{s+1}(c)^\top \hat X_s \P_s(\tau)\\[2mm]
\hline
\\[-2mm]
& 1 \end{array}
\qquad  \equiv \qquad
\begin{array}{c|c}
c &a_{c\tau}^{(s)}\\[2mm]
\hline
\\[-2mm]
& 1\end{array}~,
\end{equation}
which is clearly equivalent to (\ref{RKs1}).

\begin{rem} We observe that, in a sense, (\ref{RK1s}) can be regarded as a continuous extension of the $W$-transformation in \cite[Section IV.5]{HW2002}. Moreover, by considering, in place of (\ref{RK1s}), the following Butcher tableau,
$$\begin{array}{c|c}
c & \P_s(c)^\top X_s \P_s(\tau)\\[2mm]
\hline
\\[-2mm]
& 1 \end{array}~,$$
one obtains the continuous extension of the low-rank symplectic methods in \cite{BuBu2012}.
\end{rem}

\subsubsection{Discretization}\label{discr} We conclude this section by recalling that \cite{BIT2010,BIT2012,LIMbook2016} for the polynomial $\sigma$ defined in (\ref{RKs1})-(\ref{MFE}), one has $\sigma(h)-y(h)=O(h^{2s+1})$.\footnote{One could obtain the result also  by using the {\em symplifying assumptions} for continuous-stage RK methods \cite{TS2012,TS2014,MiBu2016}.} Moreover, by approximating the integrals
$$\int_0^1 P_j(\tau)f(\sigma(\tau h))\dd\tau$$ appearing in (\ref{MFE}) by means of a Gauss-Legendre formula 
of order $2k$, one obtains a HBVM$(k,s)$ method, which retains the order $2s$ of the approximation defined by (\ref{MFE}),
for all $k\ge s$. In particular, when $k=s$, one obtains the $s$-stage Gauss-Legendre collocation method. As a result, the Butcher tableau of a HBVM$(k,s)$ method turns out to be given by
\begin{equation}\label{hbvms}
\begin{array}{c|c}
\bfc & \I_s \P_s^\top\Omega\\[2mm]
\hline
\\[-2mm]
& \bfb^\top \end{array}  
\qquad \equiv \qquad
\begin{array}{c|c}
\bfc & \P_{s+1} \hat X_s \P_s^\top\Omega\\[2mm]
\hline
\\[-2mm]
& \bfb^\top \end{array}
\qquad =: \qquad
\begin{array}{c|c}
\bfc & A=(a_{ij})\\[2mm]
\hline
\\[-2mm]
& \bfb^\top \end{array}
\,,
\end{equation}
with $\hat X_s$ the matrix defined in (\ref{PIXs}), 
\begin{equation}\label{bc}
\bfb = \pmatrix{ccc} b_1,&\dots,&b_k\endpmatrix^\top,\qquad
\bfc = \pmatrix{ccc} c_1,&\dots,&c_k\endpmatrix^\top,
\end{equation}
the vectors containing the weights and abscissae of the quadrature, respectively,\footnote{Any quadrature is in principle allowed, provided that it is enough accurate.}
\begin{equation}\label{OmIs}
\Omega = \pmatrix{ccc} b_1\\ &\ddots\\ && b_k\endpmatrix,\qquad
\I_s = \pmatrix{ccc}
\int_0^{c_1} P_0(x)\dd x & \dots &\int_0^{c_1} P_{s-1}(x)\dd x\\
\vdots & &\vdots\\
\int_0^{c_k} P_0(x)\dd x & \dots &\int_0^{c_k} P_{s-1}(x)\dd x
\endpmatrix\in\RR^{k\times s},
\end{equation}
and
\begin{equation}\label{Pr}
\P_r = \pmatrix{ccc}
P_0(c_1) & \dots &P_{r-1}(c_1)\\
\vdots & &\vdots\\
P_0(c_k) & \dots &P_{r-1}(c_k)
\endpmatrix\in\RR^{k\times r}, \qquad r=s,s+1.
\end{equation}
In particular, from (\ref{actaus}) one obtains that the entries of matrix $A$ in (\ref{hbvms}) are given by
$$a_{ij} = a_{c_ic_j}^{(s)}, \qquad i,j=1,\dots,k.$$

\section{Second order problems}\label{sop}
Inspired by \cite{TZ2018,TSZ2019} (see also \cite{LIMbook2016}), we now consider the case of special second order problems, i.e., ODE-IVPs in the form (\ref{ivp2}).  By setting $p(t)\equiv \dot q(t)$, one then obtains the following equivalent system of first order ODEs,
\begin{equation}\label{ivp2_1}
\dot q(t) = p(t), \qquad \dot p(t) = f(q(t)), \qquad t\in[0,h], \qquad  q(0) = q_0, \, \dot q(0) = p_0\,\in\,\RR^m.
\end{equation}
HBVMs have been considered for numerically solving this problem \cite{BIT2011}. We can then consider the use of
HBVM$(\infty,s)$, too. 
To begin with, by applying same steps as above, one then obtains that (\ref{ivp2_1}) can be formally written as
\begin{eqnarray*}
\dot q(ch) &=& \P_\infty(c)^\top\otimes I_m \left[\int_0^1 \P_\infty(\tau)\otimes I_m\, p(\tau h)\,\dd\tau\right], \\
\dot p(ch) &=& \P_\infty(c)^\top\otimes I_m \left[\int_0^1 \P_\infty(\tau)\otimes I_m\, f(q(\tau h))\,\dd\tau\right], \qquad c\in[0,1].
\end{eqnarray*} 
Simplifying the expressions, integrating side by side, and imposing the initial conditions, then gives
\begin{eqnarray*}
q(ch) &=& q_0 + h\int_0^1 \I_\infty(c)^\top\P_\infty(\tau)\otimes I_m p(\tau h)\,\dd\tau, \\
p(ch) &=& p_0 + h\int_0^1\I_\infty(c)^\top\P_\infty(\tau)\otimes I_m f(q(\tau h))\,\dd\tau, \qquad c\in[0,1].
\end{eqnarray*} 
Substituting the second equation in the first one, and taking into account (\ref{IXPPI})-(\ref{X}), then gives, setting $e_1=\pmatrix{ccc} 1,&0,&\dots~\endpmatrix^\top$ and considering that $\I_\infty(c)e_1=c$,
\begin{eqnarray*}
q(ch) &=& q_0 + h\int_0^1 \I_\infty(c)^\top\P_\infty(\xi)\otimes I_m p(\xi h)\,\dd\xi\\
         &=& q_0 + h\int_0^1 \I_\infty(c)^\top\P_\infty(\xi)\otimes I_m \left[p_0 + h\int_0^1\I_\infty(\xi)^\top\P_\infty(\tau)\otimes I_m f(q(\tau h))\,\dd\tau\right]\,\dd\xi\\
         &=& q_0 + h \I_\infty(c)^\top\underbrace{\int_0^1\P_\infty(\xi)\dd\xi}_{=\,e_1}\otimes  p_0 \\ &&+\, h^2 \I_\infty(c)^\top\underbrace{\int_0^1\P_\infty(\xi)\I_\infty(\xi)^\top\dd\xi}_{=X_\infty}\int_0^1\P_\infty(\tau)\otimes I_m f(q(\tau h))\,\dd\tau \\
         &=& q_0 + ch p_0 + h^2 \int_0^1 \I_\infty(c)^\top X_\infty \P_\infty(\tau)\otimes I_m f(q(\tau h))\,\dd\tau\\
         &=& q_0 + ch p_0 + h^2 \int_0^1 \P_\infty(c)^\top X_\infty^2 \P_\infty(\tau)\otimes I_m f(q(\tau h))\,\dd\tau\\
         &=:& q_0 + ch p_0 + h^2 \int_0^1 \bar{a}_{c\tau}^{(\infty)}\otimes I_m f(q(\tau h))\,\dd\tau, \qquad c\in[0,1],
\end{eqnarray*}
where, by considering that (see (\ref{X}))
\begin{equation}\label{X2}
X_\infty^2 = \pmatrix{ccccc}
\xi_0^2-\xi_1^2 &-\xi_0\xi_1 & \xi_1\xi_2\\[2mm]
\xi_0\xi_1          &-\xi_1^2-\xi_2^2 & ~~0 & \xi_2\xi_3\\
\xi_1\xi_2          & ~~0   &-\xi_2^2-\xi_3^2 & ~~0 &\ddots\\
                          &\xi_2\xi_3 &~~0 &\ddots &\ddots  \\
                          &                &\ddots&\ddots&\ddots 
                          \endpmatrix
\end{equation}
and taking into account (\ref{PIgamma}), we have set :
\begin{eqnarray}\label{bactau}
\lefteqn{\bar{a}_{c\tau}^{(\infty)} ~=~ \I_\infty(c)^\top X_\infty\P_\infty(\tau)~\equiv~ \P_\infty(c)^\top X_\infty^2\P_\infty(\tau) ~\equiv~ \frac{1}6+\frac{\xi_1}2(P_1(c)-P_1(\tau)) ~+~}\\ 
&&     -~ \sum_{j=1}^\infty \left[(\xi_j^2+\xi_{j+1}^2)P_j(c)P_j(\tau) - \xi_j\xi_{j+1}(P_{j-1}(c)P_{j+1}(\tau)+P_{j-1}(\tau)P_{j+1}(c))\right], \qquad c,\tau\in[0,1].\nonumber
\end{eqnarray}
Moreover, by setting $q_1\equiv q(h)$ and (see (\ref{ivp2_1})) $\dot q_1\equiv p(h)$, one obtains
$$\dot q_1 = \dot q_0 + h\int_0^1 f(q(ch))\dd c,$$
and, by also considering that \quad $f(q(\tau h)) = \sum_{j\ge0} P_j(\tau)\int_0^1 P_j(\xi)f(q(\xi h))\dd\xi$,\quad $\tau\in[0,1]$,
\begin{eqnarray*}
q_1 &=& q_0 + h\int_0^1 p(ch)\dd c ~=~ q_0 + h\int_0^1\left[ \dot q_0 + h\int_0^c \dot p(\tau h)\dd\tau\right]\dd c\\
&=& q_0 + h\dot q_0 +h^2\int_0^1\int_0^c f(q(\tau h))\dd\tau\dd c \\
&=& q_0 + h\dot q_0 +h^2\int_0^1\int_0^c \left[ \sum_{j\ge0} P_j(\tau)\int_0^1 P_j(\xi)f(q(\xi h))\dd\xi\right]\dd\tau\dd c\\
&=& q_0 + h\dot q_0 +h^2\int_0^1 \left[ \sum_{j\ge0} P_j(\xi)\int_0^1\int_0^cP_j(\tau)\dd\tau\dd c\right]f(q(\xi h))\dd\xi \\
&=& q_0 + h\dot q_0 +h^2\int_0^1 \left[\P_\infty(\xi)^\top \int_0^1 \I_\infty(c)\dd c\right] f(q(\xi h))\dd\xi\\
&\equiv& q_0 + h\dot q_0 +h^2\int_0^1 \bar b_\xi f(q(\xi h))\dd\xi.
\end{eqnarray*}
Next, by taking into account (\ref{IXPPI}), one obtains:
\begin{eqnarray}\nonumber
\bar b_\xi &:=& \P_\infty(\xi)^\top \int_0^1 \I_\infty(c)\dd c ~=~ \P_\infty(\xi)^\top X_\infty^\top \int_0^1 P_\infty(c)\dd c\\
&=&\P_\infty(\xi)^\top X_\infty^\top e_1 = \xi_0 -\xi_1 P_1(\xi) ~=~ 1-\xi. \label{bxi}
\end{eqnarray}
In conclusion, we can summarize the above procedure as follows (see (\ref{bactau})):
\begin{eqnarray}\nonumber
q(ch) &=& q_0 + ch \dot q_0 + h^2 \int_0^1 \bar a_{c\tau}^{(\infty)} f(q(\tau h))\dd\tau, \qquad c\in[0,1],\\ \label{cRK2}
q_1   &=& q_0 + h\dot q_0 + h^2 \int_0^1 (1-c) f(q(ch))\dd c,\\
\dot q_1 &=& \dot q_0 + h \int_0^1 f(q(ch))\dd c.\nonumber
\end{eqnarray}
In other words, we are speaking about the application of the following ``continuous-stage Runge-Kutta-Nystr\"om (RKN, hereafter) method'' for solving problem (\ref{ivp2_1}), i.e., (\ref{ivp2}) :
\begin{equation}\label{RK2}
\begin{array}{c|c}
c & \I_\infty(c)^\top X_\infty \P_\infty(\tau)\\[2mm]
\hline
\\[-3mm]
& 1-c  \\\hline
\\[-3mm]
& 1
\end{array}
\qquad \equiv \qquad
\begin{array}{c|c}
c & \P_\infty(c)^\top X_\infty^2 \P_\infty(\tau)\\[2mm]
\hline
\\[-3mm]
& 1-c  \\\hline
\\[-3mm]
& 1\end{array}\qquad \equiv \qquad
\begin{array}{c|c}
c &\bar a_{c\tau}^{(\infty)}\\[2mm]
\hline
\\[-3mm]
& 1-c  \\\hline
\\[-3mm]
& 1\end{array}~,
\end{equation}
which provides the exact solution of the problem.

\subsection{Polynomial approximation}\label{poli2}

As done for first order problems, also in this case we can consider a polynomial approximation $\sigma\in\Pi_s$ to $q$. This is done by resorting to the same finite vectors and matrices defined in (\ref{PIs}) and (\ref{PIXs}), resulting into the following continuous-stage RKN
method:
\begin{equation}\label{RK2n}
\begin{array}{c|c}
c & \I_s(c)^\top X_s \P_s(\tau)\\[2mm]
\hline
\\[-3mm]
& 1-c  \\\hline
\\[-3mm]
& 1
\end{array}
\qquad \equiv \qquad
\begin{array}{c|c}
c & \P_{s+1}(c)^\top \hat X_s X_s \P_s(\tau)\\[2mm]
\hline
\\[-3mm]
& 1-c  \\\hline
\\[-3mm]
& 1\end{array}\qquad \equiv \qquad
\begin{array}{c|c}
c &\bar a_{c\tau}^{(s)}\\[2mm]
\hline
\\[-3mm]
& 1-c  \\\hline
\\[-3mm]
& 1\end{array}~,
\end{equation}
which defines the application of the HBVM$(\infty,s)$ method for solving (\ref{ivp2}). One has, then,
\begin{eqnarray}\nonumber
\sigma(ch) &=& q_0 + ch \dot q_0 + h^2 \int_0^1 \bar a_{c\tau}^{(s)} f(\sigma(\tau h))\dd\tau, \qquad c\in[0,1],\\ \label{cRK2s}
q_1   &=& q_0 + h\dot q_0 + h^2 \int_0^1 (1-c) f(\sigma(ch))\dd c,\\
\dot q_1 &=& \dot q_0 + h \int_0^1 f(\sigma(ch))\dd c.\nonumber
\end{eqnarray}
It is well-known \cite{LIMbook2016,BI2018} that $q_1-q(h)=\dot q_1-\dot q(h) = O(h^{2s+1})$.\footnote{Also in this case, one could derive the result through the simplifying assumptions for continuous-stage RKN methods \cite{TSZ2019}.} 

\begin{rem}\label{sge2}
We observe, however, that in order for (\ref{bxi}) to hold, one must have $s\ge 2$. Conversely, one would obtain $\bar b_\xi\equiv 1$, in place of\, $\bar b_\xi = 1-\xi$.
\end{rem}
Moreover, considering that (compare with (\ref{X2}))
\begin{equation}\label{X2s}
\hat X_s X_s = \pmatrix{ccccc}
\xi_0^2-\xi_1^2 &-\xi_0\xi_1 & \xi_1\xi_2\\[2mm]
\xi_0\xi_1          &-\xi_1^2-\xi_2^2 & ~~0 & \ddots\\
\xi_1\xi_2          & ~~0   &\ddots  & \ddots &\xi_{s-2}\xi_{s-1}\\
                          &\ddots &\ddots &-\xi_{s-2}^2-\xi_{s-1}^2 &0  \\
                          &                &\xi_{s-2}\xi_{s-1}&0&-\xi^2_{s-1}\\
                          &                &                          &\xi_{s-1}\xi_s &0
                          \endpmatrix \in \RR^{(s+1)\times s},
\end{equation}
one obtains:
\begin{eqnarray}\nonumber
\lefteqn{\bar{a}_{c\tau}^{(s)} ~=~ \I_s(c)^\top X_s\P_s(\tau)~\equiv~ \P_{s+1}(c)^\top \hat X_sX_s\P_s(\tau) ~\equiv~ \frac{1}6+\frac{\xi_1}2(P_1(c)-P_1(\tau)) ~+~}\\ \nonumber
&&     -~ \sum_{j=1}^{s-2} \left[(\xi_j^2+\xi_{j+1}^2)P_j(c)P_j(\tau) - \xi_j\xi_{j+1}(P_{j-1}(c)P_{j+1}(\tau)+P_{j-1}(\tau)P_{j+1}(c))\right] ~+~\\
&&     -~\xi_{s-1}^2P_{s-1}(c)P_{s-1}(\tau)~ -~\xi_{s-1}\xi_sP_s(c)P_{s-1}(\tau), \qquad c,\tau\in[0,1],\label{bactaus}
\end{eqnarray}
in place of (\ref{bactau}).

\subsubsection{Discretization}\label{discr2}
We conclude this section by recalling that, by approximating the integrals appearing in (\ref{cRK2s}) by means of a Gauss-Legendre formula 
of order $2k$, one obtains a HBVM$(k,s)$ method, which retains the order $2s$ of the approximation defined by (\ref{cRK2s}),
for all $k\ge s$.\footnote{In particular, when $k=s$, one obtains the RKN method induced by the $s$-stage Gauss collocation method, $s\ge 2$.} The Butcher tableau of this $k$-stage RKN method turns out to be given by:
\begin{equation}\label{hbvm2s}
\begin{array}{c|c}
\bfc & \I_s X_s\P_s^\top\Omega\\[2mm]
\hline
\\[-2mm]
& \bfb^\top\circ(1-\bfc^\top)\\[2mm]
\hline
\\[-2mm]
& \bfb^\top \end{array}  
\qquad \equiv \qquad
\begin{array}{c|c}
\bfc & \P_{s+1} \hat X_s X_s \P_s^\top\Omega\\[2mm]
\hline
\\[-2mm]
& \bfb^\top\circ(1-\bfc^\top)\\[2mm]
\hline
\\[-2mm]
& \bfb^\top \end{array}  
\qquad =: \qquad
\begin{array}{c|c}
\bfc & \bar A = (\bar a_{ij})\\[2mm]
\hline
\\[-2mm]
& \bfb^\top\circ(1-\bfc^\top)\\[2mm]
\hline
\\[-2mm]
& \bfb^\top \end{array}  
\,,
\end{equation}
with $\circ$ the Hadamard (i.e., componentwise) product, and the same matrices and vectors defined in (\ref{PIXs}) and (\ref{bc})--(\ref{Pr}). As in the case of first order problems, one has that the entries of the Butcher matrix $\bar A$ in (\ref{hbvm2s}) are given by (see (\ref{bactaus})) 
$$\bar a_{ij} = \bar a_{c_ic_j}^{(s)}, \qquad i,j=1,\dots,k,$$
for all $k\ge s$ and $s\ge 2$.

\section{Conclusions}\label{fine} In this paper, we have studied the formulation of the class of energy-conserving methods named {\em Hamiltonian Boundary Value Methods (HBVMs)} as continuous-stage RK methods. When applied for solving special second-order problems, such methods also provide a class of continuous-stage RKN methods, whose derivation has been provided in full details. The formulation of HBVMs as continuous-stage RK/RKN methods, in turn, is interesting by itself, even though the efficient implementation and analysis of the methods is better addressed, in our opinion, in their original formulation (see, e.g., the monograph \cite{LIMbook2016} or the review paper \cite{BI2018}.)

\end{document}